# The Real Non-attractive Fixed Point Conjecture and Beyond

Rajen Kumar and Tarakanta Nayak

*School of Basic Sciences, Indian Institute of Technology Bhubaneswar, India*

E-mail: rajen_2021ma04@iitp.ac.in; tnayak@iitbbs.ac.in

**Abstract.** Is it always true that every polynomial $P$ with degree at least two has a fixed point $z_0$, the real part of whose multiplier is bigger than or equal to 1, i.e., $\Re(P'(z_0)) \geq 1$? This question, raised by Coelho and Kalantari in *How many real attractive fixed points can a polynomial have?* Math. Gaz. 103 (2019), no. 556, 65–76, [3] is answered affirmatively not only for all polynomials but also for all rational functions with a superattracting fixed point. However this is not true for all rational functions. Some further investigation on distribution of multipliers of fixed points is made. Quadratic and cubic polynomials, all of whose multipliers have real part 1 are characterized. A necessary and sufficient condition is found for cubic and quartic polynomials to have all the multipliers equidistant from 1.

**2020** AMS Subject Classification: 30C10, 37F10.

## 1. Introduction

Throughout this article, $\widehat{\mathbb{C}}$ denotes the Riemann sphere $\mathbb{C} \cup \{\infty\}$ and $R : \widehat{\mathbb{C}} \to \widehat{\mathbb{C}}$ denotes a rational map, given by $R(z) = \frac{P(z)}{Q(z)}$, where $P$ and $Q$ are polynomials with complex coefficients and without any common factor. Degree of $R$ is defined as the maximum of the degrees of $P$ and $Q$. Each rational map is analytic in $\widehat{\mathbb{C}}$. It is important to note that every analytic function from $\widehat{\mathbb{C}}$ onto itself is a rational map [[2], Chapter 2]. Here the notion of analyticity at $\infty$ is as follows. If $f(\infty) = \infty$, we say $f$ is analytic at $\infty$ if $\frac{1}{f\left(\frac{1}{z}\right)}$ is analytic at 0. If $f(\infty)$ is a finite complex number, we say $f$ is analytic at $\infty$ if $f\left(\frac{1}{z}\right)$ is analytic at 0. If $f(z_0) = \infty$ for some $z_0 \in \mathbb{C}$, $f$ is said to be analytic at $z_0$ if $\frac{1}{f(z)}$ is analytic at $z_0$. A point $z_0 \in \widehat{\mathbb{C}}$ is called a fixed point of $R$ if $R(z_0) = z_0$. For $z_0 \in \mathbb{C}$, the number $\lambda = R'(z_0)$ is called the multiplier of $z_0$. If $R(\infty) = \infty$ then the multiplier of $\infty$ is defined as $h'(0)$ where $h(z) = \frac{1}{R\left(\frac{1}{z}\right)}$. The multiplier of a fixed point usually controls the local iterative behaviour of the function at the fixed point. A rational map is a polynomial if and only if it has a single pole and that is at infinity. This work is motivated by the following conjecture on the fixed points of polynomials, made by Coelho and Kalantari in [3].

**Conjecture [RNFP Conjecture].** Every polynomial with degree at least two has a fixed point, the real part of whose multiplier is bigger than or equal to 1.

The authors refer it as *the real non-attractive fixed point conjecture*. They have shown this to be true for all quadratic and cubic polynomials. We have observed that the conjecture is not only true for all higher degree polynomials but for all rational functions with a superattracting fixed point (i.e., with multiplier 0) also. For a polynomial, $\infty$ is a fixed point with multiplier 0. Therefore, the conjectured fixed point is a finite complex number for polynomials.

The assumption that the degree of the polynomial is at least two is essential. In fact, the conjecture is true for an affine map $z \mapsto az + b$ if and only if $\left|a - \frac{1}{2}\right| \leq \frac{1}{2}$ or $\Re(a) \geq 1$. To see it, note that $\infty$ is the only fixed point when $a = 1$ and its multiplier is 1. For $a \neq 1$, the fixed points are $\infty$ and $\frac{b}{1-a}$ with multipliers $\frac{1}{a}$ and $a$ respectively.

Fixed points are classified according to the value of their multipliers. A fixed point $z_0$ is called attracting if $|\lambda| < 1$ (superattracting if $\lambda = 0$). It is called indifferent if $|\lambda| = 1$. Further, it is called rationally indifferent if $\lambda^m = 1$ for some $m \in \mathbb{N}$ and irrationally indifferent if $\lambda = e^{2\pi i \alpha}$ for an irrational number $\alpha$. If $|\lambda| > 1$ then the fixed point $z_0$ is called repelling. A fixed point whose multiplier is 1 or with modulus bigger than 1 is called weakly repelling. It is known that every rational map has at least one weakly repelling fixed point (See Corollary 2.1). Weakly repelling fixed points are important in the iteration theory of rational maps. More details can be found in [1,6]. A different treatment of the subject with special focus on polynomial root finding can be found in [4]. The RNFP conjecture is about the existence of a special type of weakly repelling fixed points, namely those the real part of whose multiplier is bigger than or equal to 1.

What is needed to prove the above conjecture, is the Rational fixed point theorem. This is a well-known result. For completeness we include a detailed exposition of the latter result and a detailed discussion following Milnor [5]. Section 2 discusses the Rational fixed point theorem. Then the RNFP conjecture is shown to be true for all rational maps with at least one superattracting fixed point. In particular, this is true for all polynomials. For rational functions without any superattracting fixed point, the conjecture is sometimes, but not always true. Both types of examples are given. It also follows that every polynomial has a multiple fixed point (i.e., with multiplier 1) or has a fixed point the imaginary part of whose multiplier is non-negative. This article also deals with other issues on the distribution of multipliers. The extreme situation when all multipliers have real part equal to 1 is considered. All quadratic and cubic polynomials, all the multipliers (of fixed points) of which have real part equal to 1 are characterized in Theorem 3.1 and Theorem 3.2 respectively. A related question arises naturally: When are all the multipliers equidistant from 1? If all the fixed points of a polynomial are multiple then the answer is an obvious "yes". On the other hand if a polynomial has a simple as well as a multiple fixed point then the multipliers can never be equidistant from 1. Assuming all the fixed points to be simple Theorem 3.3 proves that if the fixed points of an $n$-degree polynomial constitute the vertices of a regular $n$-gon then all the multipliers are equidistant from 1. The converse for quadratic polynomials follows easily from the Rational fixed theorem. It is shown to be true for all cubic polynomials. For all quartic polynomials, the fixed points constitute the vertices of a rectangle whenever all the multipliers are equidistant from 1. This is Theorem 3.4. Finally, some further directions of investigation are suggested.

All rational functions considered in this article are of degree at least two unless stated otherwise. By a fixed point of a polynomial, we mean a finite fixed point of the polynomial throughout the article.

## 2. The Rational Fixed Point Theorem

A fixed point $\hat{z}$ of a rational map $R$ is a root of $R(z) - z$. The local degree $m$ of $R(z) - z$ at $\hat{z}$ is the natural number such that this map, near $\hat{z}$ behaves *like* $z \mapsto z^m$ around the origin. This number is known as the multiplicity of the fixed point $\hat{z}$.

**Definition 2.1.** *The multiplicity of a fixed point $\hat{z}$ of a rational map $R$ is defined to be the unique natural number $m$ where the Taylor series of $R(z) - z$ about $\hat{z}$ is $a_m(z - \hat{z})^m + a_{m+1}(z - \hat{z})^{m+1} + \cdots$ with $a_m \neq 0$. If $m = 1$, $\hat{z}$ is called a simple fixed point. Otherwise $\hat{z}$ is called a multiple fixed point.*

The multiplicity of a fixed point $\hat{z}$ of $R$ is the multiplicity of $\hat{z}$ considered as a root of $R(z) - z$.

The following lemma reveals the relation between the multiplier and the multiplicity of a fixed point.

**Lemma 1.** *A fixed point is multiple if and only if its multiplier is 1.*

**Proof.** Let $R(\hat{z}) = \hat{z}$ and $R(z) - z = a_m(z - \hat{z})^m + a_{m+1}(z - \hat{z})^{m+1} + \cdots$. Differentiating both sides with respect to $z$, we get $R'(z) - 1 = ma_m(z - \hat{z})^{m-1} + (m+1)a_{m+1}(z - \hat{z})^m + \cdots$. Therefore $R'(\hat{z}) = 1$ if and only if $m \geq 2$. □

Attracting, repelling and irrationally indifferent fixed points are always simple. But a rationally indifferent fixed point is simple only when the multiplier is different from 1. For example, 0 is a rationally indifferent fixed point of $iz + z^2$ and that is simple.

In order to understand all the fixed points together, it is necessary to know the number of fixed points of a rational map. The following lemma appearing in [2, Theorem 2.6.3] reveals exactly that.

**Lemma 2.** *If $R$ is a rational map of degree $d \geq 1$, then it has $d + 1$ fixed points counting multiplicity.*

**Proof.** Suppose that $R = \frac{P}{Q}$ is a rational map with degree $d \geq 1$. Then all the finite fixed points of $R$ are precisely the roots of the equation $P(z) - zQ(z) = 0$.

Suppose that the degree of $P$ is greater than that of $Q$. Then $R(\infty) = \infty$, and the degree of $P$ is $d$. Since $P(z) - zQ(z)$ is a polynomial of degree $d$, it has $d$ finite roots which are fixed points of $R$ counting multiplicity.

Suppose that the degree of $P$ is less than or equal to that of $Q$. Then $P(z) - zQ(z)$ is a poynomial of degree $d + 1$, and it has $d + 1$ finite roots counting multiplicities. Note that $\infty$ is not a fixed point of $R$ in this case. This implies that $R$ has $d + 1$ fixed points counting multiplicity. □

There is an important quantity associated with a fixed point.

**Definition 2.2.** *Given a rational map R, the residue fixed point index of R at a fixed point $\hat{z}$, denoted by $\iota(R, \hat{z})$ is defined as*

$$\frac{1}{2\pi i} \oint \frac{dz}{z - R(z)},$$

*where the integration is over a small and simple positively oriented loop around $\hat{z}$ such that it does not surround or pass through any other fixed point of R.*

The word "small and simple" in the above definition can be any circle centered at $\hat{z}$ with radius $0 < r < \min_{R(z)=z} |z - \hat{z}|$. For any such $r$, the integral has the same value by the Cauchy theorem. In fact, $\hat{z}$ is a pole of $\frac{1}{z-R(z)}$ and $\iota(R, \hat{z})$ is its residue.

Note that the analytic bijections of the Riemann sphere are precisely the Möbius maps $z \mapsto \frac{az+b}{cz+d}$, $ad - bc \neq 0$. These are important for our purpose. Two rational maps $R$ and $S$ are called conformally conjugate if there is a Möbius map $g$ such that $S(z) = g(R(g^{-1}(z)))$ for all $z \in \widehat{\mathbb{C}}$. Important normalizations are to be realized via this conjugacy. One such is used in the proof of the following lemma which gives a formula for the residue fixed point index. Even though it is given in [5], we provide a detailed proof for the convenience of the reader.

**Lemma 3.** *If the multiplier $\lambda$ of a fixed point is not equal to 1 then its residue fixed point index is given by $\frac{1}{1-\lambda}$.*

**Proof.** Without loss of generality, assume that 0 is a fixed point of a rational map $R$ with multiplier $\lambda$. If $R(\hat{z}) = \hat{z}$ and $\lambda = R'(\hat{z})$ then the map $R(z + \hat{z}) - \hat{z}$ fixes the origin and the multiplier is $\lambda$. Expanding $R$ as its Taylor series around 0, we get $R(z) = \lambda z + a_2 z^2 + a_3 z^3 + \cdots$. Since $\lambda \neq 1$, $z - R(z) = (1 - \lambda)z - a_2 z^2 - a_3 z^3 + \cdots$. This gives that

$$\frac{1}{z - R(z)} = \frac{1}{(1-\lambda)z(1 + b_1 z + b_2 z^2 + \cdots)}$$

for some complex numbers $b_1, b_2, b_3, \ldots$. Consequently, $\frac{1}{z-R(z)} = \frac{1+O(z)}{(1-\lambda)z}$ where $\frac{O(z)}{z}$ is a bounded analytic function in a neighborhood of the origin. It gives that $\iota(R, 0) = \frac{1}{2\pi i} \oint \{\frac{1}{(1-\lambda)z} + \frac{O(z)}{(1-\lambda)z}\}dz$, where the integration is over a small positively oriented circle about the origin. The residue index is nothing but $\frac{1}{1-\lambda}$ as the second term is 0 by the Cauchy theorem. □

**Remark 1.** For $\lambda = 1$, the residue index $\iota(R, \hat{z})$ is still well defined and is finite.

The following lemma proves that the residue index is invariant under conformal conjugacy. Our presentation is an elucidation of the proof already available in [5, Lemma 12.3].

**Lemma 4.** *Let $R$ and $S$ be conformally conjugate, i.e., $S(z) = g(R(g^{-1}(z)))$ for all $z \in \widehat{\mathbb{C}}$ for a Möbius map $g$. If $R(\hat{z}) = \hat{z}$ then $S(g(\hat{z})) = g(\hat{z})$ and, $\iota(R, \hat{z}) = \iota(S, g(\hat{z}))$.*

**Proof.** Note that $S(g(\hat{z})) = g(\hat{z})$ if and only if $R(\hat{z}) = \hat{z}$. Further, $R'(\hat{z}) = S'(g(\hat{z}))$. If $\hat{z}$ is a simple fixed point of $R$ then $g(\hat{z})$ is a simple fixed point of $S$ and the residue indices $\iota(R, \hat{z})$ and $\iota(S, g(\hat{z}))$ are the same. Let $\hat{z}$ be a multiple fixed point of $R$. Then $R'(\hat{z}) = 1$ and, therefore $S'(\hat{w}) = 1$ where $g(\hat{z}) = \hat{w}$. Consider a perturbed map $S_\alpha(z) = S(z) + \alpha$ for $\alpha$ sufficiently near to the origin. Also, consider the punctured disc $D_\epsilon = \{w : 0 < |w - \hat{w}| < \epsilon\}$ such that $S'(w) \neq 1$ for any $w \in D_\epsilon$, $S$ has no fixed point in $D_\epsilon$ and $S_\alpha$ has no fixed point on the boundary of $D_\epsilon$. This is possible as the zeros of a non-constant analytic function are isolated. Further, by decreasing $\epsilon$ if necessary so that $g^{-1}$ has no pole on the boundary of $D_\epsilon$, we can assume that the boundary of $g^{-1}(D_\epsilon)$ is a circle. This is because each Möbius map takes circles not passing through its pole onto circles. All the fixed points $a_1, a_2, \ldots, a_m$ of $S_\alpha$ lying in $D_\epsilon$ are simple and,

$$\sum_{i=1}^{m} \iota(S_\alpha, a_i) = \oint_{\partial D_\epsilon} \frac{dz}{z - S_\alpha(z)}.$$

Since the integral tends to $\iota(S, \hat{w})$ as $\alpha \to 0$,

$$\lim_{\alpha \to 0} \sum_{i=1}^{m} \iota(S_\alpha, a_i) = \iota(S, \hat{w}). \tag{1}$$

Let $R_\alpha(z) = g^{-1}(g(R(z)) + \alpha)$ be the perturbed map of $R$. Then $g^{-1}(a_i)$ is a simple fixed point of $R_\alpha$ with multiplier $R'_\alpha(g^{-1}(a_i))$, which is the same as $S'_\alpha(a_i)$. This gives that $\iota(R_\alpha, g^{-1}(a_i)) = \iota(S_\alpha, a_i)$ for each $i$. Consequently,

$$\sum_{i=1}^{m} \iota(S_\alpha, a_i) = \sum_{i=1}^{m} \iota(R_\alpha, g^{-1}(a_i)).$$

The right hand side is $\oint_\gamma \frac{dz}{z - R_\alpha(z)}$ where $\gamma$ is the boundary of $g^{-1}(D_\epsilon)$. Note that $\gamma$ is a circle not surrounding any fixed point of $R_\alpha$ other than $g^{-1}(a_i)$'s, and $\gamma$ does not contain any fixed point of $R_\alpha$. Further, the point $\hat{z}$ is the only fixed point of $R$ surrounded by $\gamma$. Since $R_\alpha \to R$ as $\alpha \to 0$, $\lim_{\alpha \to 0} \sum_{i=1}^{m} \iota(R_\alpha, g^{-1}(a_i)) = \iota(R, \hat{z})$. This along with Equation (1) gives that $\iota(R, \hat{z}) = \iota(S, \hat{w})$. □

The following is known as the Rational fixed point theorem, also called the Holomorphic fixed point formula. Though the proof is available, see, for example [5, Theorem 12.4], we present it in detail for the sake of completeness.

**Theorem 2.1.** *For every non-identity and non-constant rational map $R$, $\sum_{z=R(z)} \iota(R, z) = 1$.*

**Proof.** Conjugating $R$ by a Möbius map $\mu$, if necessary, we may assume without loss of generality that $R(\infty) \neq 0, \infty$. This is justified by Lemma 4. For this, $\mu$ can be chosen such that $R(\mu^{-1}(\infty)) \notin \{\mu^{-1}(0), \mu^{-1}(\infty)\}$ and then consider $\mu R \mu^{-1}$ instead of $R$. If $R(\infty) = 0$ then $\mu(z) = z + 1$ can be chosen. If $R(\infty) = \infty$ then consider $\mu(z) = \frac{1}{z-z_0}$ where $R(z_0) \neq z_0$ and $R(z_0) \neq \infty$. Such a $z_0$ can always be found for every non-constant and non-identity rational map.

Now $\lim_{z \to \infty} R(z) = l \in \mathbb{C}\setminus\{0, \infty\}$ and it means that for every $\epsilon > 0$, there is a positive number $M$ such that $|R(z) - l| < \epsilon$ whenever $|z| > M$. A bigger $M$ can be chosen, if necessary, to additionally ensure that $|R(z)| < \frac{|z|}{2}$ and all the fixed points of $R$ have modulus less than $M$. Thus it follows that $|z - R(z)| \geq |z| - |R(z)| = \frac{|z|}{2} + \frac{|z|}{2} - |R(z)| > \frac{|z|}{2}$ for all $z$ with $|z| > M$. Note that there is a $K > 0$ such that $|R(z)| < K$ for all $z, |z| > M$. Consequently, $\left|\frac{R(z)}{z(z-R(z))}\right| < \frac{2K}{|z|^2}$ whenever $|z| > M$. Let $r > M$ and $\partial D_r$ be the simple positively oriented loop $|z| = r$. We have

$$\left|\frac{1}{2\pi i}\oint_{\partial D_r} \frac{R(z)dz}{z(z-R(z))}\right| \leq \frac{1}{2\pi}\oint_{\partial D_r} \frac{|R(z)||dz|}{|z(z-R(z))|}$$
$$\leq \frac{1}{2\pi}\oint_{\partial D_r} \frac{2K|dz|}{|z|^2} = \frac{2K}{r}.$$

As $r \to \infty$, the right hand side goes to zero giving that $\lim_{r \to \infty} \frac{1}{2\pi i}\oint_{\partial D_r} \frac{R(z)dz}{z(z-R(z))} = 0$. Since $\frac{R(z)}{z(z-R(z))} = \frac{1}{z-R(z)} - \frac{1}{z}$,

$$\lim_{r \to \infty} \frac{1}{2\pi i}\oint_{\partial D_r} \frac{dz}{z-R(z)} = \lim_{r \to \infty} \frac{1}{2\pi i}\oint_{\partial D_r} \frac{dz}{z}. \quad (2)$$

Here, both the integrals are independent of $r$, for sufficiently large $r$ and the value of the right hand side is 1. Further,

$$\sum_{z=R(z)} \iota(R, z) = \sum_{j} \frac{1}{2\pi i}\oint_{\partial D_j} \frac{dz}{z-R(z)}, \quad (3)$$

where $\partial D_j$ is a small positively oriented loop that surrounds only one fixed point of $R$ and the sum is taken over all fixed points of $R$. Since all the fixed points of $R$ have moduli less than $r$, it follows from the Cauchy theorem that the integral in Equation (3) is equal to $\frac{1}{2\pi i}\oint_{\partial D_r} \frac{dz}{z-R(z)}$. Now, it follows from Equation (2) and (3) that $\sum_{z=R(z)} \iota(R, z) = 1$. □

**Remark 2.** Since the residue index of each polynomial $P$, with degree at least 2, at $\infty$ is 1, we have

$$\sum_{P(z)=z; z \in \mathbb{C}} \iota(P, z) = 0. \quad (4)$$

The following describes a relation between the type of a simple fixed point and its residue index.

**Lemma 5.** *A fixed point with multiplier $\lambda \neq 1$ is attracting if and only if the real part of its residue index is bigger than $\frac{1}{2}$.*

**Proof.** A fixed point $\hat{z}$ with multiplier $\lambda$ is attracting if $|\lambda| < 1$. Note that the function $\frac{1}{1-z}$ maps the unit disk $\{z : |z| < 1\}$ onto the half plane $H = \{z : \Re(z) > \frac{1}{2}\}$. This is clear since the function $z \mapsto 1 - z$ maps the unit disk onto the disk $D = \{z : |z - 1| < 1\}$, and the disk $D$ is mapped onto $H$ by $\frac{1}{z}$. The converse follows easily. □

**Remark 3.** From the above lemma and its proof it follows that the residue index at an indifferent fixed point has real part equal to $\frac{1}{2}$, and that of a repelling fixed point has real part less than $\frac{1}{2}$.

Here is an important consequence of the Rational fixed point theorem. This is already known [5, Corollary 12.7].

**Corollary 2.1.** *Every rational map of degree $d \geq 2$ has a repelling fixed point or a rationally indifferent fixed point with multiplier 1, or both.*

**Proof.** If possible, for a rational map $R$, suppose that there is neither a repelling fixed point nor a fixed point with multiplier 1. Since there is no fixed point of multiplier 1, all the fixed points are simple. The real part of their residue indices is bigger than or equal to $\frac{1}{2}$ by Lemma 5 and the remark following it. By the Rational fixed point Theorem, $\sum_{z=R(z)} \iota(R, z) = 1$. There are $d + 1$ simple fixed points and hence $d + 1$ terms in the sum. Consequently, the right hand side has real part bigger than or equal to $\frac{d+1}{2}$. However this is not possible since $d \geq 2$. Therefore, each rational map has a repelling fixed point or a rationally indifferent fixed point with multiplier 1, or both. □

The following theorem settles the RNFP conjecture in its full generality.

**Theorem 2.2.** *Every rational map with degree at least two and with a superattracting fixed point has a fixed point, the real part of whose multiplier is bigger than or equal to 1. In particular, this is true for all polynomials with degree at least two.*

**Proof.** If a rational map has a multiple fixed point then its multiplier is 1 and we are done. Let $R$ be a rational map, all of whose fixed points are simple. If $R$ is of degree $d$ then it has $d + 1$ distinct fixed points. Suppose that $\lambda_i = R'(z_i)$ where $R(z_i) = z_i$ for $i = 1, 2, 3, \ldots, d + 1$. By the Rational fixed point theorem, $\sum_{i=1}^{d+1} \frac{1}{1-\lambda_i} = 1$. By assumption, $R$ has a superattracting fixed point. Letting it to be $z_{d+1}$, we have $\lambda_{d+1} = 0$, and

$$\sum_{i=1}^{d} \frac{1}{1-\lambda_i} = 0. \quad (5)$$

Comparing the real part, we get that

$$\sum_{i=1}^{d} \frac{\Re(1 - \bar{\lambda}_i)}{|1-\lambda_i|^2} = 0. \quad (6)$$

Since $\Re(1 - \bar{\lambda}_i) = 1 - \Re(\lambda_i)$ and $|1 - \lambda_i|^2 > 0$ for each $i$, it is not possible to have $\Re(\lambda_i) < 1$ for all $i$. Therefore, there exists $j \in \{1, 2, \ldots, d\}$ such that $\Re(\lambda_j) \geq 1$. □

A number of important remarks follow.

**Remark 4.**

1. If $\Re(\lambda_i) > 1$ for each $i$ then each term in the left hand side of Equation (6) is negative leading to a contradiction. Therefore every polynomial with degree at least two has a fixed point, the real part of whose multiplier is less than or equal to 1. Further, if all the fixed points are simple and some fixed points have their multipliers with real part bigger than 1 then there is a fixed point whose multiplier has real part less than 1.
2. Comparing the imaginary parts in Equation (5), it is observed that $R$ has a fixed point with multiplier 1 or has a fixed point with imaginary part at least 0. In other words, if all the fixed points of a rational map are simple then it has a fixed point, the imaginary part of whose multiplier is non-negative.
3. The RNFP conjecture is true for all rational maps with a superattracting fixed point. However, nothing can be said in general for a rational map without any superattracting fixed point. For example, the two rational functions $\frac{1}{z^{d-1}}, d > 1$ and $\frac{kz}{z^2+z+1}, k \neq 0$ donot have any superattracting fixed point. Multiplier of each fixed point of the first function is $1 - d$ which is negative for $d > 2$ giving that the RNFP conjecture is false for $\frac{1}{z^{d-1}}, d > 2$. On the other hand, 0 is a fixed point of $\frac{kz}{z^2+z+1}$ with multiplier $k$. The RNFP conjecture is true for $\frac{kz}{z^2+z+1}$ whenever $\Re(k) \geq 1$.

## 3. Beyond the Real Non-attractive Fixed Point Conjecture

We now deal with some issues arising out of the real non-attractive fixed point conjecture.

### 3.1 Multipliers with real part one

When do all the finite fixed points of a polynomial have real part equal to 1? In some sense, this is an extreme situation of the RNFP conjecture. It is not difficult to construct polynomials for which the answer is yes. Let $P(z) = z + ik(z - a_1)(z - a_2) \cdots (z - a_m)(z - b_1)^{p_1}(z - b_2)^{p_2} \cdots (z - b_n)^{p_n}$ where $a_1, a_2, \ldots, a_m, b_1, b_2, \ldots, b_n, k \in \mathbb{R}\setminus\{0\}$ and $p_1, p_2, \ldots, p_n$ are natural numbers bigger than 1. Then each $a_i$ and $b_j$ are fixed points of $P$, and the real part of each multiplier is 1.

We now look at the quadratic polynomials. Since every quadratic polynomial is conformally conjugate to $z^2 + c$ for some $c$ and the multipliers are preserved under conformal conjugacy, it is enough to consider $z^2 + c$.

**Theorem 3.1.** *Let $P(z) = z^2 + c$. If $c = \frac{1}{4}$ then $P$ has a single fixed point and its multiplier is 1. If $c \neq \frac{1}{4}$ then $P$ has two simple fixed points. Further, the multipliers have real part equal to 1 if and only if $c > \frac{1}{4}$. In other words, the multiplier of each fixed point of $P$ has real part equal to 1 if and only if $c \geq 1/4$.*

**Proof.** The finite fixed points of $P(z) = z^2 + c$ are $\frac{1 \pm \sqrt{1-4c}}{2}$. There is only one fixed point if and only if $c = \frac{1}{4}$. In this case, the multiplier is 1.

If $c \neq \frac{1}{4}$ then $\frac{1+\sqrt{1-4c}}{2}$ and $\frac{1-\sqrt{1-4c}}{2}$ are two distinct fixed points with multipliers $1 + \sqrt{1 - 4c}$ and $1 - \sqrt{1 - 4c}$ respectively. The real part of these multipliers is 1 if and only if $1 - 4c < 0$ which is nothing but $c > \frac{1}{4}$. □

For dealing with general cubic polynomials, we first need to prove a simple case.

**Lemma 6.** *Suppose that a cubic polynomial $Q$ has three distinct finite fixed points $0$, $1$ and $\alpha$. Then all multipliers have real part $1$ if and only if $\alpha$ is real and $Q'(0) - 1$ is purely imaginary.*

**Proof.** If a cubic polynomial $Q$ fixes three distinct points $0$, $1$ and $\alpha$ then $Q(z) = z + kz(z-1)(z-\alpha)$ for some non-zero $k$. Further, the multipliers are $1+k\alpha$, $1+k-k\alpha$ and $1+k\alpha^2-k\alpha$. Each has real part $1$ if and only if $k\alpha$, $k - k\alpha$ and $k\alpha^2 - k\alpha$ are purely imaginary, This is true if and only if $k$ is purely imaginary and $\alpha$ is real. The proof completes by observing that $k = \frac{Q'(0)-1}{\alpha}$. □

**Theorem 3.2.** *Suppose that a cubic polynomial has three distinct fixed points. Then the real part of all the multipliers is $1$ if and only if all the fixed points are collinear and the real part of any one multiplier is $1$.*

**Proof.** If $a, b$ and $c$ are the distinct fixed points of a cubic polynomial $P$ then $Q = \phi P \phi^{-1}$ fixes $0$ and $1$ where $\phi(z) = \frac{z-a}{b-a}$. The other fixed point of $Q$ is $\frac{c-a}{b-a}$, which is clearly different from $0$ and $1$. Note that $P'(a) = Q'(0)$, $P'(b) = Q'(1)$ and $P'(c) = Q'\left(\frac{c-a}{b-a}\right)$. Further, $a, b, c$ are collinear if and only if $0, 1, \frac{c-a}{b-a}$ are so. These follow from the facts that the inverse of an affine map is affine, and an affine map takes lines onto lines. Letting $\alpha = \frac{c-a}{b-a}$, observe that $Q(z) = z + kz(z-1)(z-\alpha)$ for some non-zero $k$.

If all the multipliers (of fixed points) of $P$ have real part equal to $1$ then the same is true for $Q$. Since the multipliers of $Q$ are $Q'(0) = 1+k\alpha$, $Q'(1) = 1+k-k\alpha$ and $Q'(\alpha) = 1 + k\alpha^2 - k\alpha$, $k\alpha$ and $k - k\alpha$ are purely imaginary. This happens only when $k$ is purely imaginary and $\alpha$ is real. In other words, the arguments of $c - a$ and $b - a$ are same or differ by $\pi$. Therefore $a, b, c$ are collinear. Obviously, the real part of a (in fact every) multiplier is $1$.

Conversely, let the fixed points $a, b, c$ of $P$ be collinear and one multiplier has real part equal to $1$. Let the fixed point of $P$ whose multiplier has real part $1$ be denoted by $a$. Then the fixed points $0, 1, \alpha = \frac{c-a}{b-a}$ of $Q$ are collinear, giving that $\frac{c-a}{b-a}$ is real. Since the real part of $P'(a)$ is $1$, the real part of $Q'(0) = 1 + k\alpha$ is $1$. In other words, $k$ is purely imaginary. It is non-zero since $Q'(0) \neq 1$. Now, by Lemma 6, the real part of each multiplier of $Q$ is $1$. Therefore, all multipliers of $P$ have real part equal to $1$. □

**Remark 5.** Up to conformal conjugacy, the only cubic polynomial (with only simple fixed points) whose all multipliers have real part $1$ is $kz^3 - (k + k\alpha)z^2 + (k\alpha + 1)z$ for some non-zero and purely imaginary $k$ and a real number $\alpha$ different from $0$ and $1$.

### 3.2 Multipliers equidistant from 1

Equation (6) in Section 2 gives that $\sum_{i=1}^{d} \mathfrak{Re}(1 - \bar{\lambda}_i) = 0$ whenever all the multipliers are equidistant from $1$. In this situation, the average of the real parts of all multipliers is $1$. We now look at the situation when all the multipliers are equidistant from $1$. If all the fixed points of a polynomial are multiple then every multiplier is $1$. Such a polynomial is $P(z) = z + (z - a_1)^{p_1}(z - a_2)^{p_2} \cdots (z - a_k)^{p_k}$ for some $a_1, a_2, \ldots, a_k \in \mathbb{C}$ and positive integers $p_1, p_2, \ldots, p_k$, each of which is bigger than $1$. In this case, all the finite fixed points $a_1, a_2, \ldots, a_k$ have the same multiplier, namely $1$, and all the multipliers are trivially equidistant from $1$. However, for polynomials with a multiple as well as a simple fixed point, all multipliers can never be equidistant from $1$. Thus, the question of all multipliers being equidistant from $1$ makes sense only when all the fixed points are simple.

The case is again straightforward for quadratic polynomials. If every fixed point of a quadratic polynomial is simple then both multipliers, $\lambda_1, \lambda_2$ are equidistant from $1$. In fact, it follows from the Rational fixed point theorem that $\lambda_1 + \lambda_2 = 2$, which gives $|\lambda_1 - 1| = |\lambda_2 - 1|$. Polynomials with degree at least three are considered in the next theorem. For stating this we need a definition. For a natural number $n$, a regular $n$-gon is a (convex) polygon with vertices at $v_k = a + re^{i(\theta + \frac{2\pi k}{n})}$, $k = 0, 1, 2, \ldots, (n-1)$ for some $a \in \mathbb{C}, r > 0$ and $\theta \in (0, 2\pi]$. A regular $n$-gon is completely determined by $a, r$ and $\theta$. We say it is centered at $a$ and with radius $r$.

**Theorem 3.3.** *Let $P$ be a polynomial with degree $n$, $n \geq 3$ and each of its fixed points is simple. If all the fixed points form the vertices of a regular $n$-gon then all the multipliers are equidistant from $1$. Moreover, its multipliers constitute the vertices of a regular $n$-gon centered at $1$.*

**Proof.** If all the fixed points of a polynomial $P$ with degree $n$ are simple, and form the vertices of a regular $n$-gon centered at $a$ and with radius $r$ then $P(z) = z + M((z-a)^n - (re^{i\theta})^n)$ for some non-zero complex number $M$. Note that each $v_k = a + re^{i(\theta + \frac{2\pi k}{n})}$ is a fixed point of $P$ and $P'(v_k) = 1 + Mn(re^{i(\theta + \frac{2\pi k}{n})})^{n-1}$, for $k = 0, 1, 2, \ldots, n-1$ which is nothing but $1 + Mn(re^{i\theta})^{n-1} e^{-\frac{2\pi ik}{n}}$. Note that each $P'(v_k)$ is equidistant from 1. Further, these multipliers constitute the vertices of a regular $n$-gon with centre at 1 and radius $|M|nr^{n-1}$. □

It is natural to ask when the converse of the above result is true. For cubic and quartic polynomials, we answer this question.

**Theorem 3.4.** *Let $P$ be a polynomial with degree $n$ and each of its fixed points is simple. If all the multipliers are equidistant from* 1 *then all the fixed points are on the vertices of an equilateral triangle or a rectangle when $n = 3$ or $4$ respectively.*

**Proof.** Let all the fixed points of a polynomial with degree $n$ be simple. Then

$$P(z) = z + k(z - \alpha_1)(z - \alpha_2)(z - \alpha_3) \cdots (z - \alpha_n), \quad (7)$$

for some distinct complex numbers $\alpha_1, \alpha_2, \ldots, \alpha_n$ and a non-zero complex number $k$. Differentiating (7), we get

$$P'(z) = 1 + k \sum_{i=1}^{n} L_i(z) \quad \text{where} \quad L_i(z) = \prod_{j=1, j \neq i}^{n} (z - \alpha_j). \quad (8)$$

Note that $L_i(\alpha_j) = 0$ for $i \neq j$ and $\lambda_i = P'(\alpha_i) = 1 + kL_i(\alpha_i)$. Now, $|1 - \lambda_i| = |1 - \lambda_j|$ if and only if $|L_i(\alpha_i)| = |L_j(\alpha_j)|$ for each $i \neq j$. Note that $|L_i(\alpha_i)|$ is nothing but the product of distances of all other fixed points from $\alpha_i$.

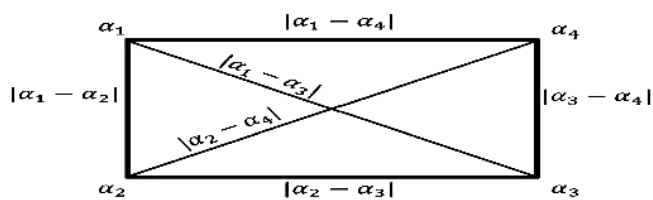

Figure 1. Rectangle with vertices $\alpha_1, \alpha_2, \alpha_3$ and $\alpha_4$

If $n = 3$ then the assumption of this theorem means that $|(\alpha_1 - \alpha_2)(\alpha_1 - \alpha_3)| = |(\alpha_2 - \alpha_1)(\alpha_2 - \alpha_3)| = |(\alpha_3 - \alpha_1)(\alpha_3 - \alpha_2)|$. In other words, $|\alpha_1 - \alpha_2| = |\alpha_2 - \alpha_3| = |\alpha_3 - \alpha_1|$. Hence the fixed points form the vertices of an equilateral triangle.

Let $n = 4$ and $\alpha_1, \alpha_2, \alpha_3$ and $\alpha_4$ be the fixed points of $P$. Then, denoting the distance $|\alpha_i - \alpha_j|$ by $\alpha_{ij}$, it is observed that $\alpha_{12}\alpha_{13}\alpha_{14} = \alpha_{21}\alpha_{23}\alpha_{24}$, $\alpha_{21}\alpha_{23}\alpha_{24} = \alpha_{31}\alpha_{32}\alpha_{34}$, $\alpha_{31}\alpha_{32}\alpha_{34} = \alpha_{41}\alpha_{42}\alpha_{43}$ and $\alpha_{41}\alpha_{42}\alpha_{43} = \alpha_{12}\alpha_{13}\alpha_{14}$. The first two equations give that $\alpha_{12}\alpha_{14} = \alpha_{32}\alpha_{34}$ whereas it follows from the second and third that $\alpha_{21}\alpha_{23} = \alpha_{41}\alpha_{43}$. Since $\alpha_{ij} > 0$ and $\alpha_{ij} = \alpha_{ji}$, these equalities give that $\alpha_{14} = \alpha_{23}$ and consequently $\alpha_{12} = \alpha_{34}$. Now, the fixed points $\alpha_i, i = 1, 2, 3, 4$ can be thought of as the vertices of a quadrilateral, whose opposite sides have same length. In other words, this is a parallelogram. Further, it also follows from $\alpha_{12}\alpha_{13}\alpha_{14} = \alpha_{21}\alpha_{23}\alpha_{24}$ that $\alpha_{13} = \alpha_{24}$, i.e., the length of diagonals of this parallelogram are same giving that it is a rectangle. □

We conclude by stating some possible directions of research.

For a polynomial with at least two multiple fixed points, the corresponding multipliers are 1. It is natural to ask when two simple fixed points have the same multiplier. Even though never true for any quadratic polynomial, this is possible for certain higher degree polynomials. For example, the polynomial $z^4 - 2z^2$ fixes 0 and $-1$, each of whose multiplier is 0. However a general condition ensuring this may be found.

**Question 1.** When the multipliers of two (or all) simple fixed points of a polynomial are same?

Theorem 2.2 and Remark (3) following it show that some rational maps satisfy the RNFP conjecture whereas some do not. In view of this, the following seems worth doing.

**Question 2.** Characterize all the rational maps for which the RNFP conjecture is true.

The RNFP conjecture can also be made for periodic points. A point $z_0$ is called a $p$-periodic point of a rational function $R$, if $p$ is the smallest natural number such that $z_0$ is a fixed point of $R^p$. But a fixed point of $R^p$ may not be a $p$-periodic point of $R$. It follows from Theorem 2.2 that for $p > 1$, there is a fixed point of $P^p$ whose multiplier has real part at least 1. However, it does not seem to be trivial to decide whether this is a $p$-periodic point of $P$, and this leads to:

**Question 3.** Given polynomial $P$ and $p > 1$, does there exist a $p$-periodic point of $P$ whose multiplier has real part equal to or bigger than 1?

Dynamics of polynomials, stated in Remark 5 may be studied. For every polynomial $P$, there is an open set $A$ such that $P^n(z) \to \infty$ as $n \to \infty$ for all $z \in A$. The boundary of this set is known as the Julia set of $P$. The approximate Julia sets of $kz^3 - (k + k\alpha)z^2 + (k\alpha + 1)z$ are given in Figure 2a and 2b for $k = 1.2 = \frac{1}{\alpha}$ and $k = 0.01 = \frac{1}{\alpha}$ respectively. It is the boundary of the black region. In both the figures, the Julia set is disconnected whereas it is totally disconnected only in the second figure. This calls for further investigation.

**Question 4.** Study the Julia sets of polynomials all of whose fixed points have multipliers with real part 1.

The same question can be asked for polynomials satisfying the assumptions of Theorem 3.3.

**Question 5.** Study the Julia sets of polynomials all the multipliers of whose fixed points are equidistant from 1.

As noted earlier, the RNFP conjecture is true for all rational maps having a superattracting fixed point. But the local degree of the rational map $R$ at this fixed point may not be equal to the degree of $R$. This is where $R$ differs from a polynomial, and the iterative behaviour (the Fatou and the Julia set) of such rational maps calls for further study. The Fatou set of a rational map $R$ is the set of all points in a neighbourhood of which the sequence of iterates $\{R^n\}_{n>0}$ is equicontinuous [2, Chapter 2]. The Julia set of $R$ is the complement of its Fatou set in $\widehat{\mathbb{C}}$. The following can be undertaken.

**Question 6.** Study the Fatou and the Julia set of rational maps which satisfy the RNFP conjecture but are not conformally conjugate to any polynomial.

## Acknowledgement

The second author is supported by a MATRICS project (MTR/2018/000498 dated 12.03.2019), funded by SERB, Govt of India.

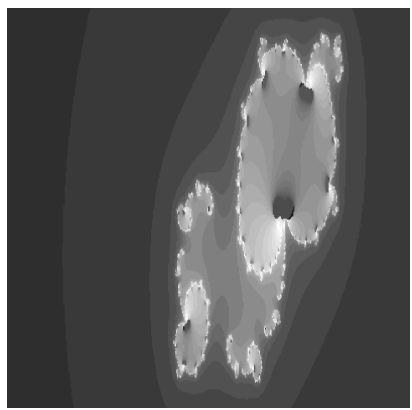

Figure 2a. The Julia set of $1.2i(z^3 - z^2) + i(z - z^2) + z$

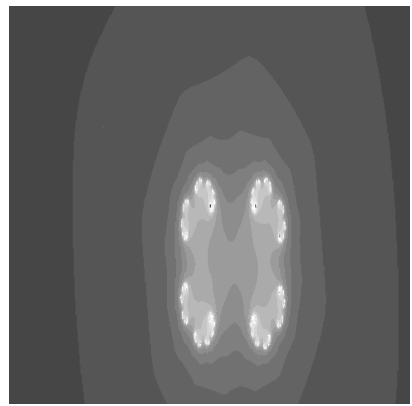

Figure 2b. The Julia set of $0.01i(z^3 - z^2) + i(z - z^2) + z$